# Hyperelliptic Curves over Small Finite Fields and GPU Accelerators

Martin Raum[1]

**Abstract:** We present a hardware-accelerated computation of Hasse-Weil invariants of all hyperelliptic curves of given genus over a fixed finite field. Our main motivation is the determination of traces of Frobenius on cohomology corresponding moduli stacks à la Bergström-Faber-van-der-Geer. This paper also constitutes a case study of the performance of lookup table based Zech arithmetic on Graphics Processing Units (GPUs). GPUs are by now ubiquitous in numerics, to an extend that their development itself is propelled by scientific computing. Algebraic computing has profited very little from the advancement of hardware design. We suggest that specific computations can be benificially adjusted to GPUs with comparatively little effort.

**hyperelliptic curves** ■ **Hasse-Weil invariants** ■ **finite fields** ■ **Zech logarithms** ■ **hardware accelerators** ■ **GPUs**
**MSC Primary: 14Q05** ■ **MSC Secondary: 11G25, 14H45**

IN a series of papers Bergström, Faber, and van der Geer showed how to obtain traces of Frobenius on cohomology of local systems on moduli stacks of curves and of abelian varieties of genus 2 and 3 from (weighted) multiplicities of Hasse-Weil invariants on the moduli stack of genus 2 and 3 curves over finite fields [Ber09; BFG08; BFG14; BG08; FG04a; FG04b]. From their computations emerged presumable motivic pieces in the cohomology of $\mathcal{M}_{3,n}$ that are not associated with automorphic forms for $\mathrm{Sp}_3$. Instead, a series of subsequent computations by Chenevier, Renard, Lannes, Taïbi, and Mégarbané [CL18; CR15; Még16; Taï17] suggests that they do match automorphic forms for $\mathrm{SO}_7$ and $\mathrm{SO}_9$. However, it remains unexplained why. More experimental data might give inspiration towards a resolution of this mystery.

Counts of Hasse-Weil invariants of hyperelliptic curves over finite fields $\mathbb{F}_q$ constitute the computational heart of Bergström's, Faber's, and van der Geer's findings. The first motivation for the present work is to extend the range of previous data in genus 2 and 3 and to include the case of genus 4. We make this possible by a more detailed understanding of how finite field arithmetics performs on hardware accelerators. As a result, we can provide data[2] for odd prime powers $q$, if $q < 300$ and $g = 2$, if $q < 50$ and $g = 3$, and if $q < 15$ and $g = 4$. Data that we generated shall eventually be incorporated into a database of Siegel modular forms [BFG17] in the cases of $g = 2$ and $g = 3$. The high point of our computation is the case of $g = 4$, which bears the potential to detect new motivic pieces in the cohomology of moduli $\mathcal{H}_{4,n}$ of marked hyperelliptic curves by means of [Ber09]. To compute traces on cohomology of $\mathcal{M}_{3,n}$, we will have to also consider plane curves of degree 4 [Ber08]; cohomology of $\mathcal{M}_{4,n}$ requires even further considerations [Tom05]. Both cases will be pursued separately in a sequel.

---

The author was partially supported by Vetenskapsrådet Grant 2015-04139.

[2]This computation was partially performed using resources of the Chalmers Centre for Computational Science and Engineering C3SE.

By work of Harvey, Hasse-Weil invariants can be computed in asymptotically polynomial time [Har14; HS14; HS16], and a more than excellent implementation of the resulting algorithm is available in Sutherland's smalljac. Harvey's algorithm, however, requires a certain precomputation that must be performed for each curve. As a result, it is not practical for the small finite fields over which our hyperelliptic curves are defined. The naive determination of Hasse-Weil invariants by evaluation of a Weierstraß model at all points turns out to be superior. Performance can then be improved by tweaking arithmetics of finite fields.

This connects to our second motivation for this paper: We explore the potential towards hardware accelerated arithmetics in finite fields. We provide an implementation of our central Algorithm 1 for both CPUs and GPUs in order to compare them. It is based on lookup tables for Zech logarithms. On our test system we are able to verify a speedup by a factor of 17 in the relevant parameter range when using GPUs; and 37 for larger primes. Beyond that we improve by several orders of magnitude on the performance of the previously employed computer program.

In our discussion of performance, we evaluate the impact of two issues: Random memory access patterns arising from lookup tables and branch divergence arising from a dichotomy in our data types. While the latter does not influence performance severely, the former sets the limits of our approach. In light of the rapid improvements of memory and cache that is built into scientific GPUs, this is good news as it suggests that the relative performance advantage of our approach against mere CPU-focused implementations might further increase in the close future.

Finite fields are crucial for many constructions in theoretical mathematics. While, for example, algebraic geometry over finite fields features them intrinsically, they appear as auxiliary objects in multi-modular algorithms. Characteristics of finite fields that appear in the context of theoretical mathematics are rather different from those arising in elliptic curve cryptography or in RSA crypto-systems. In particular, their sizes are moderate, as opposed to excessively large fields in cryptography that aim at making the computation of discrete logarithms virtually impossible. The medium range of field size that make appearance in our computation has apparently not yet been addressed in the context of co-processors. In this paper, we close this gap. Given that revised programming standards like OpenMP 4.5 make the utilization of hardware accelerators easier than ever before, we expect that our work might be instrumental to a broader adoption within, for example, the overlap of computational mathematics with number theory and various flavors of geometry.

*Acknowledgement*   The author is grateful to Jonas Bergström for helpful discussions and his comments on an early version of this manuscript.

# 1 Counts of hyperelliptic curves and Hasse-Weil invariants

We summarise the information on hyperelliptic curves that we need to formulate Algorithm 1. Compare [Ber09] for a similar discussion. We throughout fix a power $q$ of an odd prime $p$ and a genus $g \in \mathbb{Z}_{\geq 2}$.

**1.1 Hasse-Weil invariants of hyperelliptic curves**  Fix a finite field $\mathbb{F}_q$ with $q$ elements. A smooth and proper algebraic curve $C$ of genus $g$ is called hyperelliptic if it admits a morphism to $\mathbb{P}^1$ of degree 2. Since $q$ is odd, there is a defining equation

$$C : Y^2 = f(X) = \sum_{n=0}^{2g+2} c_n X^n$$

for a polynomial $f(X)$ of degree $2g+1$ or $2g+2$. After homogenizing $f$, this yields a description of $C$ as a subset of the weighted projective space $\mathbb{P}(1, g+1, 1)$ including possible points $(1 : y : 0)$ at infinity. If $f$ is square-free then $C$ is smooth. For convenience, we set $f(\infty) := c_{2g+2}$.

Given a square-free polynomial $f(X)$ over $\mathbb{F}_q$ with associated hyperelliptic curve $C_f$ and a finite field $\mathbb{E} \supseteq \mathbb{F}_q$, we have

$$\#C_f(\mathbb{E}) = r^0(f, \mathbb{E}) + 2r^\square(f, \mathbb{E}) + \bigl(1 + \chi_\mathbb{E}(f(\infty))\bigr),$$

where

$$r^0(f, \mathbb{E}) := \#\{x \in \mathbb{E} : f(x) = 0\},$$
$$r^\square(f, \mathbb{E}) := \#\{x \in \mathbb{E} : f(x) \text{ is a square in } \mathbb{E}\},$$

and $\chi_\mathbb{E}$ is the quadratic character associated with $\mathbb{E}$. We write

$$a(f, \mathbb{E}) = 1 + \#\mathbb{E} - \#C_f(\mathbb{E}) \tag{1.1}$$

for the Hasse-Weil invariants of $C_f$.

From the decomposition of $\#C_f(\mathbb{E})$ into $r^0(f, \mathbb{E})$ and $r^\square(f, \mathbb{E})$, it becomes clear that

$$\#C_f(\mathbb{E}) + \#C_{f^-}(\mathbb{E}) = 2 + \#\mathbb{E} \quad \text{for} \quad f^-(X) = f(-X).$$

We call $f^-$ the quadratic twist of $f$, following the established terminology for the associated curves $C_f$ and $C_{f^-} = (C_f)^-$.

Let $\mathrm{Pg}_g(\mathbb{F}_q)$ be the set of square-free polynomials of degree $2g+1$ and $2g+2$ with coefficients in $\mathbb{F}_q$. We write $\Lambda(d)$ for the set of integer partitions of $d \in \mathbb{Z}_{\geq 0}$. Given $f \in \mathrm{Pg}_g(\mathbb{F}_q)$ of degree $2g+2$, we let $\lambda(f) = \deg f_1 + \cdots + \deg f_l$ be the integer partition associated with the factorization $f(X) = \prod_{i=1}^l f_i$ into irreducible polynomials over $\mathbb{F}_q$. If $\deg f = 2g+1$, we let $\lambda(f)$ be that partition amended by an additional summand 1.

The goal of this paper is to compute for fixed $g$ and $q$ the following map:

$$\begin{aligned}
\mathrm{N}_{g,q} : \mathbb{Z}_{\geq 0}^g \times \Lambda(2g+2) &\longrightarrow \mathbb{Z}_{\geq 0}, \\
\bigl((a_1, \ldots, a_g), \lambda\bigr) &\longrightarrow \#\{f \in \mathrm{Pg}_g(\mathbb{F}_q) : \forall 1 \leq e \leq g : a(f, \mathbb{F}_{q^e}) = a_e, \lambda(f) = \lambda\}.
\end{aligned} \tag{1.2}$$

Bergström describes how to extract from its values the traces of Frobenius on cohomology of the moduli stack of marked hyperelliptic curves [Ber09].

**1.2 Affine reduction of hyperelliptic curves**  There is an action of $\mathrm{GL}_2(\mathbb{F}_q) \times \mathbb{F}_q^\times$ on $\mathrm{Pg}_g(\mathbb{F}_q)$ that preserves $\#C_f(\mathbb{E})$:

$$\left(f\big|\big(\left(\begin{smallmatrix}a & b\\ c & d\end{smallmatrix}\right), e\big)\right)(X) := (cX+d)^{2g+2} f\Big(\frac{aX+b}{cX+d}\Big)/e^2.$$

Its orbits are in general hard to determine, but the action of $G(\mathbb{F}_q) = \mathrm{Aff}_1(\mathbb{F}_q) \times \mathbb{F}_q^\times$ is transparent. Here

$$\mathrm{Aff}_1(\mathbb{F}_q) := \left\{\left(\begin{smallmatrix}a & b\\ 0 & 1\end{smallmatrix}\right) \in \mathrm{GL}_2(\mathbb{F}_q)\right\}$$

is the group of affine transformations. We next discuss a reduction theory for this action. To this end, fix sets of representatives $R_i$ of $\mathbb{F}_q^\times / \mathbb{F}_q^{\times i}$, $i \in \mathbb{Z}_{\geq 0}$.

Fix $d \in \{2g+1, 2g+2\}$ and let $P_d(\mathbb{F}_q)$ be the set of square-free polynomials of degree $d$ with coefficients in $\mathbb{F}_q$. Recall that $q$ is a power of a prime $p$. We consider the cases $p \nmid d$ and $p \mid d$ separately. Suppose that $p \nmid d$ and consider a polynomial $f(X) \in P_d(\mathbb{F}_q)$. Replacing $X$ by $X - c_{d-1}/dc_d$, we see that every orbit of $G(\mathbb{F}_q) \circlearrowright P_d(\mathbb{F}_q)$ contains at least one representative $f(X) = \sum c_n X^n$ such that $c_{d-1} = 0$. This reduction corresponds the action of $\left(\left(\begin{smallmatrix}1 & b\\ 0 & 1\end{smallmatrix}\right), 1\right) \in G(\mathbb{F}_q)$.

Set

$$i := \inf\{0 \leq i' < d : c_{i'} \neq 0\}, \quad j := \inf\{0 \leq j' < i : c_{j'} \neq 0\}.$$

If $j \neq -\infty$, using the action of $\left(\left(\begin{smallmatrix}a & 0\\ 0 & 1\end{smallmatrix}\right), 1\right) \in G(\mathbb{F}_q)$, we may replace $X$ by $aX$ and thus achieve that $c_j/c_i \in R_{i-j}$. We have $i \neq -\infty$, since $f$ is square-free and $d \geq 5$. After employing the action of $\left(\left(\begin{smallmatrix}1 & 0\\ 0 & 1\end{smallmatrix}\right), e\right) \in G(\mathbb{F}_q)$, we may therefore assume that $c_i \in R_2$. This shows that if $p \nmid d$, then $P_d(\mathbb{F}_q)$ equals $G(\mathbb{F}_q)\widetilde{P}_d(\mathbb{F}_q)$ where $\widetilde{P}_d(\mathbb{F}_q)$ is the intersection of $P_d(\mathbb{F}_q)$ with

$$\big\{c_d X^d + c_0 X^0 : c_d \in \mathbb{F}_q^\times, c_0 \in R_2\big\}$$
$$\cup \bigcup_{\substack{1 \leq i \leq d-2 \\ 0 \leq j < i}} \left\{\begin{array}{l} c_d X^d + c_i X^i + c_j X^j + \cdots + c_0 X^0 : \\ c_d \in \mathbb{F}_q^\times, c_i \in R_2, c_j \in c_i R_{i-j}, c_{j-1}, \ldots, c_0 \in \mathbb{F}_q \end{array}\right\}.$$

We shall refer to the elements of $\widetilde{P}_d(\mathbb{F}_q)$ as reduced. The number of polynomials $P_d(\mathbb{F}_q)$ with reduction $f(X) = \sum c_i X^i$ will be called the multiplicity of $f$ and is denoted by $\mathrm{mult}(f)$. It clearly only depends on the support of the coefficients. Specifically, we have

$$\mathrm{mult}(c_d X^d + c_0 X^0, \mathbb{F}_q) = \frac{\#\mathbb{F}_q \#\mathbb{F}_q^\times}{\#R_2} = \frac{q(q-1)}{2},$$

$$\mathrm{mult}(c_d X^d + c_i X^i + c_j X^j + \cdots + c_0 X^0, \mathbb{F}_q) = \frac{\#\mathbb{F}_q (\#\mathbb{F}_q^\times)^2}{2\#R_{i-j}} = \frac{q(q-1)^2}{2\gcd(q-1, i-j)}.$$

Now consider the case $p \mid d$ and let $f \in P_d(\mathbb{F}_q)$. If $c_{d-1} \neq 0$, then we can replace $X$ by $X - c_{d-2}/(d-1)c_{d-1}$ to ensure that $c_{d-2} = 0$. The remaining reduction process is the same as before, so that we may choose as a set of reduced polynomials $\widetilde{P}_d(\mathbb{F}_q)$ the square-free ones in the following set:

$$\left\{c_d X^d + c_{d-1} X^{d-1} + c_0 X^0 : c_d \in \mathbb{F}_q^\times, c_{d-1} \in \mathbb{F}_q, c_0 \in R_2\right\}$$
$$\cup \bigcup_{\substack{1 \le i \le d-3 \\ 0 \le j < i}} \left\{\begin{array}{c} c_d X^d + c_{d-1} X^{d-1} + c_i X^i + c_j X^j + \cdots + c_0 X^0 : \\ c_d \in \mathbb{F}_q^\times, c_{d-1} \in \mathbb{F}_q, c_i \in R_2, c_j \in c_i R_{i-j}, c_{j-1}, \ldots, c_0 \in \mathbb{F}_q \end{array}\right\}.$$

The corresponding numbers of elements in $P_d(\mathbb{F}_q)$ with prescribed reduction are

$$\mathrm{mult}(c_d X^d + c_{d-1} X^{d-1}, \mathbb{F}_q) = q,$$
$$\mathrm{mult}(c_d X^d + c_{d-1} X^{d-1} + c_0 X^0, \mathbb{F}_q) = \frac{q(q-1)}{2},$$
$$\mathrm{mult}(c_d X^d + c_{d-1} X^{d-1} + c_i X^i + c_j X^j + \cdots + c_0 X^0, \mathbb{F}_q) = \frac{q(q-1)^2}{2\gcd(q-1, i-j)}.$$

**1.3 Quadratic twists**  We have already noted that Hasse-Weil invariants of quadratic twists can be computed from one another. Since $P_d(\mathbb{F}_q)$ consists of polynomials with at least one coefficient that varies in $R_2 = \mathbb{F}_q^\times / \mathbb{F}_q^{\times 2}$, we can easily make use of this observation by replacing $\widetilde{P}_d(\mathbb{F}_q)$ with $\widetilde{P}'_d(\mathbb{F}_q)$: the intersection of $P_d(\mathbb{F}_q)$ with the following sets:

$$\left\{c_d X^d + X^0 : c_d \in \mathbb{F}_q^\times\right\}$$
$$\cup \bigcup_{\substack{1 \le i \le d-2 \\ 0 \le j < i}} \left\{\begin{array}{c} c_d X^d + X^i + c_j X^j + \cdots + c_0 X^0 : \\ c_d \in \mathbb{F}_q^\times, c_j \in c_i R_{i-j}, c_{j-1}, \ldots, c_0 \in \mathbb{F}_q \end{array}\right\}$$

if $p \nmid d$, and otherwise

$$\left\{c_d X^d + c_{d-1} X^{d-1} + X^0 : c_d \in \mathbb{F}_q^\times, c_{d-1} \in \mathbb{F}_q\right\}$$
$$\cup \bigcup_{\substack{1 \le i \le d-3 \\ 0 \le j < i}} \left\{\begin{array}{c} c_d X^d + c_{d-1} X^{d-1} + X^i + c_j X^j + \cdots + c_0 X^0 : \\ c_d \in \mathbb{F}_q^\times, c_{d-1} \in \mathbb{F}_q, c_j \in c_i R_{i-j}, c_{j-1}, \ldots, c_0 \in \mathbb{F}_q \end{array}\right\}.$$

**1.4 The algorithm**  We suggest Algorithm 1 to compute $N_{g,q}$. It is structured with an eye to our map-reduce based implementation. Our intention is to perform line 7–14 on a GPU. We elaborate on details in Section 3.

# 2 Arithmetics of finite fields

Line 8 of Algorithm 1 requires us to evaluate a polynomials over finite field. To efficiently do so, we have employed Zech based arithmetics on GPU-accelerators. In this section, we revisit the concept of Zech logarithms and some of the state-of-the-art implementations of finite field arithmetic.

---

**Algorithm 1:** Counting hyperelliptic curves
**Data:** odd prime power $q$, genus $g \geq 1$
**Result:** curve counts $N_{q,g}$

1  $N_{g,d}(\bullet, \bullet) \leftarrow 0$;
2  **for** $d \in \{2g+1, 2g+2\}$ **do**
3     **for** $c_d X^d + \cdots + c_0 X^0 \in \widetilde{P}'_d(\mathbb{F}_q)$ **do**
4        $r^0_{1 \leq \bullet \leq g} \leftarrow 0$;
5        $r^\square_{1 \leq \bullet \leq g} \leftarrow 0$;
6        **for** $1 \leq e \leq g$ **do**
7           **for** $x \in \mathbb{F}_{q^e}$ **do**
8              $f \leftarrow \sum_{i=0}^{d} c_i x^i$;
9              **if** $f = 0$ **then**
10                $r^0_e \leftarrow r^0_e + 1$;
11             **else if** $f \in \mathbb{F}^2_{q^e}$ **then**
12                $r^\square_e \leftarrow r^\square_e + 2$;
13             **end**
14          **end**
15       **end**
16       **for** $1 \leq e \leq g$ **do**
17          compute $a_e \leftarrow a_e(\sum c_i X^i)$ from $r^0_e$ and $r^\square_e$;
18       **end**
19       compute $\lambda \leftarrow \lambda(\sum c_i X^i)$ from $r^0_\bullet$ if possible; otherwise by factoring;
20       $m \leftarrow \mathrm{mult}_q(\sum c_i X^i)$;
21       increment $N_{g,d}((a_1, \ldots, a_g), \lambda)$ by $m$;
22       increment $N_{g,d}((1+q-a_1, \ldots, 1+q^g-a_g), \lambda)$ by $m$;
23    **end**
24 **end**

---

**2.1 Finite fields of prime order and machine integers** We start by finite fields $\mathbb{F}_p$ of prime order whose elements can be straightforwardly represented by machine integers provided that $p$ is not too large. Any $a \in \mathbb{F}_p$ has a unique preimage $0 \leq \tilde{a} < p$ under the projection $\mathbb{Z} \twoheadrightarrow \mathbb{F}_p$. Such a representation of $a$ as an integer is called normalized. On modern architectures this requires $p < 2^{32}$ or $p < 2^{64}$, i.e.

$$p \leq 4294967291 \quad \text{or} \quad p \leq 18446744073709551557$$

*Addition and subtraction* To accommodate addition, it is common to assume that $2p < 2^{32}$ or $2p < p^{64}$. In this way, addition can be performed by three instructions: (1) Addition $\tilde{a} + \tilde{b}$, (2) Comparison $\tilde{a} + \tilde{b} \geq p$. (3) Subtraction $\tilde{a} + \tilde{b} - p$ if normalization is necessary.

If $4p < 2^{32}$ or $4p < 2^{64}$ this can be tuned by switching the order of comparison and subtraction. Further improvements are possibly by delaying normalization. Addition in finite fields implemented in this way can be freely assumed to be fast.

Subtraction can be handle by a similar scheme and turns out to be slightly faster than addition, as it does not require normalization.

*Multiplication*   Multiplication in finite fields gives rise to two challenges. First, multiplication of normalized representatives in general exceeds the capacity of 32- and 64-bit integers. This is not too serious an issue, since efficient representation in a pair of CPU registers is possible. Second, reduction modulo $p$ cannot be substituted by a bare comparison as in the case of addition. Since the computation of integer remainders is costly, this has major bifurcations on how multiplication is implemented. Based on precomputed inverses one can significantly lower the impact, but nonetheless multiplication consumes a multitude of instruction-cycles.

The situation is generally worse on GPUs, which are even less optimized for the computation of integer remainders than CPUs.

**2.2 General finite fields and machine integers**   Elements of $\mathbb{F}_q$, $q = p^e$ cannot be represented naively by integers if $e > 1$, but $\mathbb{F}_q$ can be represented as a quotient of a polynomial ring with indeterminate $T$:

$$\mathbb{F}_q = \mathbb{F}_p[T] / R_q(T) \mathbb{F}_p[T]$$

for a suitable polynomial $R_q(T)$ of degree $e$. Finding a good choice of $R_q(T)$ is not straightforward, but we will not be concerned with this problem. In order to optimize multiplication, it is desirable to employ $R_q(T)$ most of whose coefficients vanish.

In analogy with the case of $\mathbb{F}_p$, we can represent elements of $\mathbb{F}_q$ by polynomials

$$c_0 + c_1 T^1 + \cdots c_{e-1} T^{e-1}, \quad 0 \leq c_i \leq p-1.$$

Addition and subtraction thus reduce to the corresponding operation on $e$ elements of $\mathbb{F}_p$. This is compatible with bit-packing, that increases memory efficiency and thus performance.

Multiplication requires us to take the product of polynomials and subsequently reduce them. Needless to say, both are inherently slow operations.

**2.3 Finite fields and Zech representations**   The group of units $\mathbb{F}_q^\times = \mathbb{F}_q \setminus \{0\}$ in $\mathbb{F}_q$ is a cyclic group, i.e. there is a generator $g_q \in \mathbb{F}_q$ such that every element $b \in \mathbb{F}_q^\times$ can be represented as $b = g_q^k$ for some $0 \leq k \leq q-2$. This observation is at the foundation of Zech representations for finite fields: Nonzero elements are encoded as exponents $k$ for one fixed choice of $g_q$.

As for multiplication the Zech representation allows for a straightforward reduction to addition modulo $q-1$, except in the case of multiplication by zero:

$$g_q^k \cdot g_q^l = g_q^{k+l} \quad \text{and} \quad g_q^k \cdot 0 = 0.$$

Addition in Zech representations is reduced to one multiplication (hence, to one addition modulo $q-1$) and the computation of so-called Zech logarithms:

$$g_q^k + g_q^l = g_q^k \cdot (1 + g_q^{k-l}) \quad \text{and} \quad g_q^k + 0 = g_q^k \cdot 1 = g_q^k \cdot g_q^0. \tag{2.1}$$

Given $g_q$ and $0 \le k \le q-2$ the Zech logarithm of $k$ is the unique $0 \le m(k) \le q-2$ such that

$$g_q^{m(k)} = 1 + g_q^k,$$

provided that the right hand side is nonzero in $\mathbb{F}_q$.

It should be clear that Zech logarithms are expensive to determine, as they manifest discrete logarithms. For this reason, Zech logarithms are typically precomputed and are exclusively employed for small $q$. Lookup tables are of size $\mathcal{O}(q)$, improving on $\mathcal{O}(q^2)$ for lookup tables of multiplication in more traditional representations of elements of $\mathbb{F}_q$.

**2.4 Implementations of arithmetic in finite field** It seems impossible to exhaustively discuss available implementations of finite field arithmetic. We restrict our discussion to CPU-based libraries FLINT [FLINT], NTL [NTL], and FFLAS [FFLAS], GPU-based studies in [GIT09; LMA12; MP10; TYS14], and [DSC14] for a discussion of implementations on FPGAs.

FLINT is a library for number theory written in C. It provides three implementations of finite fields, based on data types `fq`, `fq_nmod`, and `fq_zech`. The first two represent elements of $\mathbb{F}_q = \mathbb{F}_p[T]/R_q(T)$ as polynomials in $T$ with coefficients in $\mathbb{F}_p$. In the first case, coefficients are stored in an arbitrary precision format, while the second one makes the assumption that $p$ be small enough to fit them into machine size integers. Addition and multiplication is realized by arithmetic of polynomials over $\mathbb{Z}$ with subsequent reduction modulo $p$ and $R_q(T)$, if necessary. The third implementation available in FLINT is based on Zech representations, and relies on precomputed Zech logarithms.

NTL is a library for number theory written in C++. It provides types `ZZ_pEX` and `ZZ_pEX` serving the same purpose as `fq` and `fq_nmod` in FLINT. Representations are essentially the same, but use of different algorithms and implementations impacts runtime performance, which Victor Shoup summarized in a recent note[3]. Most notatbly, NTL supports thread parallelism for some of the computations, e.g. for computing in $\mathbb{F}_q$ as opposed to $\mathbb{F}_p$.

FFLAS resulted from an effort to utilize floating point arithmetic for the purpose of computing in finite fields $\mathbb{F}_p$ of prime size. Bounds on rounding errors allow to execute several arithmetic operations in a row without correcting the numeric representation of corresponding integers. The library is implemented in C++ and uses extensive templating to provide flexibility concerning the represenations; The Givaro library provides alternative, more classical representations of finite field elements.

One of the early evaluations of GPU-based prime field arithmetic is contained in [GIT09]. The authors focused on primes $p$ of 200 to 600 bits, targeting elliptic curve cryptography. They could realize a speedup of about 2.6, after addressing problems of how to keep integers of corresponding size in the GPU registers.

Little later [MP10] investigated the case of polynomial multiplication over prime fields on GPUs. While they argue that computations over general finite fields can be reduced to the prime field case, one can also view their contribution as applicable to general finite fields when using a representation in terms of polynomials. The authors summarize that they achieved speedups of about 30.

---

[3]<http://www.shoup.net/ntl/benchmarks.pdf>

Continuing in some sense the efforts of [GIT09], multiplication from a more hardware-centered perspective was investigated in [LMA12]. Tying individual GPU threads more closely to each other, they arrived at speedups between 27.7 and 71.4 as a result.

Work in [TYS14] is most closely related to the present paper. The authors considered the evaluation of quadratic, multivariate polynomials over fields $\mathbb{F}_q$ with $q = 2^e$. However, their choice of $e = 32$ is prohibitively large for lookup tables and Zech representations.

Recall that the standard implementation of multiplication in finite fields mostly suffers from the slow determination of integer remainders. A custom design on FPGAs can mitigate this problem. A concrete description of circuits reflecting standard algorithms was given in [DSC14].

## 3 Implementation

We have provided an implementation HyCu[4] of Algorithm 1 in C++ and OpenCL. In this section we describe some of its details. It should be noted that HyCu dates back to 2016, a time prior to the release and adoption in standard compilers of OpenMP 4.5. The perspective that features of OpenMP 4.5 facilitate enormously the use of techniques described here, contributed to the author's motivation to disseminate the experiences he has made.

**3.1 Data types** Recall that we have fixed a prime power $q$. We use the same encoding of finite field elements in the CPU and GPU code. We assume throughout that $2q$ does not exceed machine word size. As in Section 2.3, we can fix a generator $g_q$ of $\mathbb{F}_q$. In practice, this choice is implicitly made via the finite field interface of FLINT. We encode elements $b \in \mathbb{F}_q$ by

$$b \longmapsto \begin{cases} q-1, & \text{if } b = 0; \\ k, & \text{if } b \neq 0 \text{ and } b = g_q^k. \end{cases} \tag{3.1}$$

**3.2 Arithmetic via lookup tables** Recall from Section 2.3 that Zech based arithmetic requires the computation of remainders modulo $q-1$ and of Zech logarithms. We provide both via lookup tables `exp_red_table` and `incr_table`, respectively. The encoding in (3.1) requires one special case. Specifically, we tabulate in an array `incr_table` the following function:

$$k \longmapsto \begin{cases} q-1, & \text{if } 1 + a_q^k = 0; \\ m(k), & \text{if } 1 + a_q^k = a_q^{m(k)}; \end{cases} \quad \text{for } 0 \leq k \leq q-2. \tag{3.2}$$

Multiplication of numbers $a$ and $b$ can be implemented in a straightforward way as

```
if ( a == prime_power_pred || b == prime_power_pred )
  c = prime_power_pred;
else
  c = exp_red_table[a + b];
```

---

[4] https://github.com/martinra/hycu

Here `prime_power_pred` stores $q-1$. Notice that we face one possible branch divergence, when running this on a GPU. Section 4 contains a brief discussion.

Addition is realized along the lines of

```
if ( a == prime_power_pred )        c = b;
else if ( b == prime_power_pred )   c = a;
else if ( a > b ) {
  tmp = incr_table[a − b];
  c = tmp == prime_power_pred ? prime_power_pred :
                                exp_red_table[tmp + b];
} else {
  tmp = incr_table[b − a];
  c = tmp == prime_power_pred ? prime_power_pred :
                                exp_red_table[tmp + a];
}
```

In light of (2.1) it is natural to ask whether tabulating $m(k)$ for $2-q \le k \le q-2$ (as opposed to $0 \le k \le q-2$) could speedup addition by saving one comparison and thus avoiding further branch divergence. We have experimented with such a change and to our own surprise could not detect any speedup.

**3.3 Array parallelization and queuing of polynomials**   Our implementation is naturally divided into CPU and GPU hosted code. The amount of GPU specific code is relatively small. The loops in lines 2, 3, and 6 of Algorithm 1 are run on the CPU. We do not discuss them further. The computation of $r_e^0$ and $r_e^\square$ in lines 7–14 is carried out on a GPU via OpenCL. The computation is done via a map-reduce algorithm. Write $r_e^0(x)$ and $r_e^\square(x)$ for the contributions of $x$ to $r_e^0$ and $r_e^\square$ in lines 10 and 12. Then we perform the following transformations:

$$[x : x \in \mathbb{F}_{q^e}^\times] \rightsquigarrow [(r_e^0(x), r_e^\square(x)) : x \in \mathbb{F}_{q^e}^\times] \rightsquigarrow \left(r_e^0 = \sum_{x \in \mathbb{F}_{q^e}^\times} r_e^0(x), r_e^\square = \sum_{x \in \mathbb{F}_{q^e}^\times} r_e^\square(x)\right). \tag{3.3}$$

The first transformation is array parallelized on the GPU. Its output is stored in two arrays corresponding to $r^0(x)$ and $r^\square(x)$, whose indices correspond to $x$ via the representation in (3.1). For our later discussion of performance, it is important to record that the associated OpenCL kernel exclusively uses global arrays.

```
evaluate(
  global const int * restrict poly_coeffs_exp,
  const int prime_power_pred,
  global const int * restrict exp_red_table,
  global const int * restrict incr_table,
  global int * restrict nmbs_unramified,
  global int * restrict nmbs_ramified
  )
```

The value of $x$ is obtained as an array index

```
int x = get_global_id(0);
```

In particular, we exploit the fact that our encoding of $x \in \mathbb{F}_q^\times$ allows for a natural enumeration via array indices between 0 and $q-2$. We exclude the case $x = 0$ from the GPU-based computation to avoid one possible branch divergence in the computation.

The second transformation in (3.3) is based on the GPU reduction implementation from a 2010 white paper[5]. In the range of $p$ that we are interested in, reduction consumes significantly less time than evaluation. We therefore confine ourselves to a remark that naive GPU based reduction can be sped up significantly by carrying out the last few reduction steps on the CPU.

## 4 Performance

HyCu contains the GPU based implementation of $N_{d,q}$ from (1.2) and a CPU based alternative. To evaluate and compare performance, however, we use the Hasse-Weil invariants (1.1) of a single curve. They are also available in smalljac[6], which is optimized for the case of large $q$, and therefore give a better impression of performance relative to alternative approaches. Our test system features one NVidia Tesla K80 as a GPU and two Intel Xeon E5-2683 as CPUs. We will use only one core of the latter.

Table 1 contains minimal runtimes for the computation of $a(f, \mathbb{E})$ over various finite fields for a single curve defined over $\mathbb{Z}$. We give the time needed to compute the lookup tables from Section 3.2 in the second column. They are used in both the CPU and the GPU implementation, which are detailed in the third through fifth column. Timings for our GPU implementation are split into runtimes for the evaluation and the reduction step described in Section 3.3. The next to last column of Table 1 displays runtimes for a naive implementation based on FLINT's Zech representations of finite field elements. The final column provides timings for smalljac.

*Outline*  Runtime of the evaluation kernel contains an apparent outliner at $p = 1511$. However, we could consistently reproduce the measurement on our test system, and record it for completeness. Beyond that, we could observe a similar runtime for $p = 1499$ and $p = 1493$, while performance was as expected for $p = 1489$. Given that we heavily rely on lookup tables, this degraded performance can be most likely explained by cache associativity. We also observe that GPU runtime drops twice at $p = 401$ and $p = 2003$. Again, due to our predominant use of lookup tables, we assume this is connected to cache hierarchy, but could not completely explain our observation.

*Comparison with smalljac*  We see that for a single curve precomputation of Zech logarithms dominates the runtime. In the intended setting however, that contribution must be divided by $\#\widetilde{\mathrm{P}}_d(\mathbb{F}_q) \approx q^{d-2}$. Already in the case of $q = 101$, this reduces the average contribution of precomputing Zech logarithms to 45 picoseconds per curve. Smalljac, on the other hand, for small $p$ is dominated by precomputations which have to be performed for each curve, and therefore do not average out. In the range $p \le 500$ that is of greatest interest to us, we obtain an acceleration by a factor of about 17 when comparing our CPU and GPU based implementations. This factor increases to 37 at $p = 3001$, but at this point smalljac (on a single CPU thread) already outperforms our GPU implementation by a factor of more than 2. One might argue that smalljac was optimized in detail, while this is very much not the case for HyCu. One nevertheless has to keep in mind the asymptotically better runtime of smalljac.

---

[5]http://developer.amd.com/resources/articles-whitepapers/opencl-optimization-case-study-simple-reductions/
[6]http://math.mit.edu/~drew/smalljac_v4.1.3.tar

|         | tables | CPU  | GPU  |        | naive Zech | smalljac |
|---------|--------|------|------|--------|------------|----------|
| $p = q$ |        |      | eval | reduce |            |          |
| 101     | 4.70   | 0.63 | 0.36 | 0.10   | 3.04       | 43       |
| 211     | 23.8   | 4.47 | 1.47 | 0.15   | 15.9       | 43       |
| 307     | 54.9   | 12.1 | 3.11 | 0.22   | 30.9       | 43       |
| 401     | 99.2   | 22.2 | 0.95 | 0.35   | 54.0       | 43       |
| 503     | 167    | 36.3 | 1.25 | 0.83   | 87.7       | 43       |
| 1009    | 892    | 148  | 2.93 | 6.49   | 348        | 44       |
| 1511    | 2545   | 349  | 394  | 19.7   | 778        | 44       |
| 2003    | 5093   | 789  | 11.5 | 36.1   | 1462       | 44       |
| 3001    | 13474  | 3669 | 15.4 | 82.6   | 5091       | 44       |

Table 1: Counting points over $\mathbb{F}_{p^2}$ of $Y^2 = 12X^6 + 13X^5 + 10X^4 + 11X^3 + 16X^2 + 12X + 17$, runtime in milliseconds

*Branch divergence*   In all GPU based implementations branch divergence is a possibly worry. Representation (3.1) suggests that branch divergence could also negatively impact the performance of our implementation. Every addition in the sum $\sum c_i x^i$, $x \in \mathbb{F}_{q^e}$ can yield such a branch divergence if either $c_i = 0$ or the partial sum $\sum_{i' < i} c_{i'} x^{i'}$ vanishes. Since we iterate over all polynomials in $\widetilde{P}_d(\mathbb{F}_q)$, for a heuristic estimate, it is reasonable to assume that these cases occur with probability $1/q^e$ each. As a result, we expect that no more than a proportion of $2d/q^e$ additions yields branch divergence. Branch divergence is therefore neglectable.

*Cache behavior*   More seriously, our significant use of lookup table brings up concerns about caching behavior and memory access times. Memory characteristics of GPUs are less well understood than those for CPUs. An attempt to acquire systematic knowledge was pursued in [And+14]. Our benchmarks were carried out on a Tesla K80, for which (configurable) L1 cache is maximal 112KB. This, for example, explains the performance drop between $p = 101$ and $p = 211$. Using 16-bit tables instead of 32-bit tables could then improve performance up $p = 181$ if $g = 2$. Bit packing would possibly extend this range further.

[And+14]   M. Andersch, J. Lucas, M. Alvarez-Mesa, and B. Juurlink. "Analyzing GPGPU Pipeline Latency". *Proc. 10th Int. Summer School on Advanced Computer Architecture and Compilation for High-Performance and Embedded Systems, Fiuggi, Italy (ACACES' 14)*. 2014.

[Ber08]    J. Bergström. "Cohomology of moduli spaces of curves of genus three via point counts". *J. Reine Angew. Math.* 622 (2008).

[Ber09]    J. Bergström. "Equivariant counts of points of the moduli spaces of pointed hyperelliptic curves". *Doc. Math.* 14 (2009).

[BFG08]    J. Bergström, C. Faber, and G. van der Geer. "Siegel modular forms of genus 2 and level 2: cohomological computations and conjectures". *Int. Math. Res. Not. IMRN* (2008).

[BFG14]    J. Bergström, C. Faber, and G. van der Geer. "Siegel modular forms of degree three and the cohomology of local systems". *Selecta Math. (N.S.)* 20.1 (2014).

[BFG17]    J. Bergström, C. Faber, and G. van der Geer. *Siegel Modular Forms of Degree Two and Three*. http://smf.compositio.nl. 2017.

[BG08]     J. Bergström and G. van der Geer. "The Euler characteristic of local systems on the moduli of curves and abelian varieties of genus three". *J. Topol.* 1.3 (2008).

[CL18]     G. Chenevier and J. Lannes. *Automorphic forms and even unimodular lattices*. 2018.

[CR15]     G. Chenevier and D. Renard. "Level one algebraic cusp forms of classical groups of small rank". *Mem. Amer. Math. Soc.* 237.1121 (2015).

[DSC14]    J.-P. Deschamps, G. D. Sutter, and E. Cant. *Guide to FPGA Implementation of Arithmetic Functions*. Springer Publishing Company, Incorporated, 2014.

[FFLAS]    T. F.-F. group. *FFLAS-FFPACK: Finite Field Linear Algebra Subroutines / Package*. v2.2.2. http://linalg.org/projects/fflas-ffpack. 2016.

[FG04a]    C. Faber and G. van der Geer. "Sur la cohomologie des systèmes locaux sur les espaces de modules des courbes de genre 2 et des surfaces abéliennes. I". *C. R. Math. Acad. Sci. Paris* 338.5 (2004).

[FG04b]    C. Faber and G. van der Geer. "Sur la cohomologie des systèmes locaux sur les espaces de modules des courbes de genre 2 et des surfaces abéliennes. II". *C. R. Math. Acad. Sci. Paris* 338.6 (2004).

[FLINT]    W. Hart, F. Johansson, and S. Pancratz. *FLINT: Fast Library for Number Theory*. v2.5.2. http://flintlib.org. 2015.

[GIT09]    P. Giorgi, T. Izard, and A. Tisserand. "Comparison of modular arithmetic algorithms on GPUs". *ParCo'09: International Conference on Parallel Computing*. 2009.

[Har14]    D. Harvey. "Counting points on hyperelliptic curves in average polynomial time". *Ann. of Math. (2)* 179.2 (2014).

[HS14]     D. Harvey and A. V. Sutherland. "Computing Hasse-Witt matrices of hyperelliptic curves in average polynomial time". *LMS J. Comput. Math.* 17.suppl. A (2014).

[HS16]     D. Harvey and A. V. Sutherland. "Computing Hasse-Witt matrices of hyperelliptic curves in average polynomial time, II". *Frobenius distributions: Lang-Trotter and Sato-Tate conjectures*. Vol. 663. Contemp. Math. Amer. Math. Soc., Providence, RI, 2016.

[LMA12]    K. Leboeuf, R. Muscedere, and M. Ahmadi. "High performance prime field multiplication for GPU". *IEEE International Symposium on Circuits and Systems (ISCAS)*. 2012.

[Még16]    T. Mégarbané. *Traces des opérateurs de Hecke sur les espaces de formes automorphes de* $SO_7$, $SO_8$ *ou* $SO_9$ *en niveau* 1 *et poids arbitraire*. arXiv:1604.01914. 2016.

[MP10]   M. M. Maza and W. Pan. "Fast polynomial multiplication on a GPU". *Journal of Physics: Conference Series*. Vol. 256. 1. 2010.

[NTL]    V. Shoup. *NTL: A Library for doing Number Theory*. v10.3.0. <http://www.shoup.net/ntl/>. 2016.

[Taï17]  O. Taïbi. "Dimensions of spaces of level one automorphic forms for split classical groups using the trace formula". *Ann. Sci. Éc. Norm. Supér. (4)* 50.2 (2017).

[Tom05]  O. Tommasi. "Rational cohomology of the moduli space of genus 4 curves". *Compos. Math.* 141.2 (2005).

[TYS14]  S. Tanaka, T. Yasuda, and K. Sakurai. "Fast evaluation of multivariate quadratic polynoimals over gf (232) using graphics processing units". *Journal of Internet Services and Informations Security (JISIS)* 4.3 (2014).

Chalmers tekniska högskola och Göteborgs Universitet, Institutionen för Matematiska vetenskaper, SE-412 96 Göteborg, Sweden

E-mail: <martin@raum-brothers.eu>

Homepage: <http://raum-brothers.eu/martin>